# The factorization method for the Askey-Wilson polynomials.

Gaspard Bangerezako [*]

Institut de Mathématique, Université Catholique de Louvain,

Chemin du Cyclotron 2, B-1348 Louvain-La-Neuve, Belgium.

email: bangerezako@agel.ucl.ac.be

April 24, 1998

**Abstract**. A special Infeld-Hull factorization is given for the Askey-Wilson second order $q$-difference operator. It is then shown how to deduce a generalization of the corresponding Askey-Wilson polynomials.

## 1. Introduction.

A given second order $q$-difference operator

$$H(z;n) = u_2(qz;n)\mathbf{E}_q^2 + u_1(qz;n)\mathbf{E}_q + u_0(qz;n) \qquad (1)$$

where

$$\mathbf{E}_q^i(\Omega(z;n)) = \Omega(zq^i;n), \quad z \in R, \ i,n \in Z, \qquad (2)$$

is said to be *factorizable* according to *Infeld-Hull* (this is a $q$-version) if and only if the following product can be performed [7], [15]:

$$\begin{aligned} H(z;n) - \mu(n) &= (\eta(z;n)\mathbf{E}_q + g(z;n))\,(\zeta(z;n)\mathbf{E}_q + f(z;n)) \\ H(z;n+1) - \mu(n) &= (\zeta(z;n)\mathbf{E}_q + f(z;n))\,(\eta(z;n)\mathbf{E}_q + g(z;n)) \end{aligned} \qquad (3)$$

The operators $\eta(z;n)\mathbf{E}_q+g(z;n)$ and $\zeta(z;n)\mathbf{E}_q+f(z;n)$ are called "*lowering*" and "*raising*" respectively. It is thus a particular case of the well known Darboux transformation (where $n$ is no more a variable but all simply an index) [5], [14], [18], [21].

---

[*]*Permanent address* : Université du Burundi, Faculté des sciences, Département de Mathématique, B. P. 2700 Bujumbura, Burundi, East Central Africa.

Supposing that $u_2$ doesn't depends on $n$, let us reduce (3) in the following:

$$\hat{H}(z;n) - \hat{\mu}(n) = (u_2(z)\mathbf{E}_q + \hat{g}(z;n))(u_2(z)\mathbf{E}_q + \hat{f}(z;n))$$
$$\hat{H}(z;n+1) - \hat{\mu}(n) = (u_2(z)\mathbf{E}_q + \hat{f}(z;n))(u_2(z)\mathbf{E}_q + \hat{g}(z;n)), \quad (4)$$

where $\hat{H}(z;n) = u_2(z)H(z;n)$. The corresponding system then reads:

$$\hat{f}(qz;n+1) + \hat{g}(z;n+1) = \hat{f}(z;n) + \hat{g}(qz;n)$$
$$\hat{f}(z;n+1)\hat{g}(z;n+1) = \hat{f}(z;n)\hat{g}(z;n) + \hat{\mu}(n) - \hat{\mu}(n+1) \quad (5)$$

As a general property, the system (5) (more precisely, a system equivalent to (5)) is known to be a discrete version of the Toda lattice (see for example [23]). It is easy to see that the intertwining relations in (4) (isospectral deformations) will then constitute a discrete version of the corresponding Lax-pair. The system (5) can also be considered as a discrete delay-Painlevé I [6]. It is morever known to admit symmetries [20]. We are however in this work concerned in applications of the IHF rather than in study of its intrinsic structure.

An alternative to the IHF is the *Inui* factorization [8] where the "upstair" (raising) and "down-stair" (lowering) operators are not obtained from a direct factorization of the studied second order difference operator itself but from some difference or contiguous relations satisfied by its eigenfunctions given in a general a priori known form. An Inui type factorization for the Askey-Wilson operator [1] (see [10] for the canonical form) can be found in [16], [9] and [4] where interesting applications were also obtained. In [22] and [23] IHF techniques were applied to the three-term recurrence relations for orthogonal polynomials.

In [19] and [3] were given special IHF of the following difference operator

$$\tilde{\sigma}(s)\Delta\nabla + \tilde{\tau}(s;l)\Delta - \tilde{\lambda}(n;l) \quad (6)$$

solutions of which are the $l$-th difference of the discrete hypergeometric polynomials on a linear lattice $P_n^{(l)}(s) := \Delta^l(P_n(s))$ [17]. Here $\Delta = \mathbf{T} - 1$, $\nabla = 1 - \mathbf{T}^{-1}$ ($\mathbf{T}^i p(s) := p(s+i)$), $\tilde{\sigma}$ and $\tilde{\tau}$ being polynomials (in $s$) of degree $\leq 2$ and 1 respectively, $\tilde{\lambda}$ a constant (in $s$). In [19] the role of the "variable of factorization" (as $n$ in (3)) is played by the order $l$ of the corresponding difference ($n$ remaining fixed) while in [3] that role is played by the degree of the corresponding polynomials ($l$ remaining fixed (to zero, for simplicity)). Evidently, the IHF from [19] can be extended to the case of non-linear lattices. In that case, one intersects (formally) the evoked above Inui type factorizations from [16], [9], [4].

The main objective of the present work is to give an extension of [3] to the case of non-linear lattices. But here, to be concrete, we will work with the Askey-Wilson second order $q$-difference operator in its canonical form [10], instead of the Nikiforov-Suslov-Uvarov hypergeometric difference operator [17]. A simple and nonsignificant reduction of the lattice and an elementary change of variables lead indeed at the latter to the former. The IHF that we are going to carry out for the Askey-Wilson operator will furnish various formulas such as *strict* (i.e. without any perturbation of parameters) difference relations, recurrence relations and Rodrigues type formulas for the general sequences of eigenfunctions of the Askey-Wilson operator. Next, a somewhat converse reasoning will suffice to attain the most ambitious (see the conclusion) aim of this work: to implement the IHF system (5) as an independent source of the Askey-Wilson Polynomials [1]. All those discussions will be held in the following section. In the last one, we will extract a Rodrigues type formula allowing to produce (starting at the difference hypergeometric functions obtained in [2]) a sequence of functions generalizing the Askey-Wilson polynomials.

## 2. The factorization method for the Askey-Wilson second order $q$-difference operator.

The Askey-Wilson polynomials (AWP) [1], $\mathcal{P}_n(\chi(z))$, $\chi(z) = \frac{z+z^{-1}}{2}$ are defined by:

$$\mathcal{P}_n(\chi(z)) = \frac{(ab, ac, ad; q)_n}{a^n} {}_4\phi_3\left( \begin{array}{c} q^{-n}, abcdq^{n-1}, az, az^{-1} \\ ab, ac, ad \end{array} \bigg| q; q \right) \qquad (7)$$

where the basic hypergeometric (or $q$-hypergeometric) series ${}_r\phi_s$ read:

$${}_r\phi_s\left( \begin{array}{c} a_1, \ldots, a_r \\ b_1, \ldots, b_s \end{array} \bigg| q; z \right) := \sum_{k=0}^{\infty} \frac{(a_1, \ldots, a_r; q)_k}{(b_1, \ldots, b_s; q)_k} (-1)^{(1+s-r)k} q^{(1+s-r)\binom{k}{2}} \frac{z^k}{(q; q)_k}$$

---
[1] The original idea of generating discrete orthogonal polynomials from the IHF system (5) comes from [15] but the approach adopted there (as in [7]) already doesn't allow to handle totally the Hahn case.

with
$$(a_1,\ldots,a_r;q)_k := (a_1;q)_k \ldots (a_r;q)_k,$$
while
$$(\sigma;q)_0 := 1, (\sigma;q)_k := \prod_{i=0}^{k-1}(1-\sigma q^i), k = 1,2,3,\ldots \tag{8}$$

They satisfy the following second order $q$-difference equation [10]:
$$\mathcal{L}\mathcal{P}_n(\chi(z)) = \lambda(n)\mathcal{P}_n(\chi(z)) \tag{9}$$
where
$$\mathcal{L}(z) = v(z)\mathbf{E}_q - (v(z) + v(z^{-1})) + v(z^{-1})\mathbf{E}_q^{-1} \tag{10}$$
with
$$v(z) = \frac{(1-az)(1-bz)(1-cz)(1-dz)}{(1-z^2)(1-qz^2)}$$
$$\mathbf{E}_q^i(\mathcal{P}_n(\chi(z))) = \mathcal{P}_n(\chi(q^i z)), i \in Z, \tag{11}$$
is the Askey-Wilson second order $q$-difference operator
and
$$\lambda(n) = -(1-q^{-n})(1-abcdq^{n-1}).$$

In anticipation on the coming $q$-integration and basing ourself on the experience from [3], we are led to reformulate, before the factorization, the equation (9) as follows:
$$\{\mathcal{A}(qz)\mathbf{E}_q^2 - [\mathcal{A}(qz) + \mathcal{B}(qz) + \mathcal{K}(qz)\lambda(n)]\mathbf{E}_q + \mathcal{B}(qz)\}\mathcal{P}_n(\chi(z)) = 0 \tag{12}$$
where
$$\mathcal{A}(z) = \frac{A_{-2}z^{-2} + A_{-1}z^{-1} + A_0 + A_1 z + A_2 z^2}{qz - z^{-1}}$$
$$\mathcal{B}(z) = \frac{A_2 z^{-2} + A_1 z^{-1} + A_0 + A_{-1} z + A_{-2} z^2}{z - qz^{-1}}$$
$$\mathcal{K}(z) = z - z^{-1}$$
$$A_{-2} = 1; A_{-1} = -(a+b+c+d); A_0 = ab+ac+ad+bc+bd+cd$$
$$A_1 = -(abc+abd+bcd+acd); A_2 = abcd$$

Consider next the following $q$-difference operator
$$\mathcal{H}(z;n) = \mathcal{A}(z)\left[\mathcal{A}(qz)\mathbf{E}_q^2 - [\mathcal{A}(qz) + \mathcal{B}(qz) + \mathcal{K}(qz)\lambda(n)]\mathbf{E}_q + \mathcal{B}(qz)\right]. \tag{13}$$

The central result of this work is the following

**Proposition 1** *The operator $\mathcal{H}$ in (13) factorizes into*

$$\mathcal{H}(z;n) - \mu(n) = (\mathcal{A}(z)\mathbf{E}_q + \mathcal{G}(z;n))(\mathcal{A}(z)\mathbf{E}_q + \mathcal{F}(z;n))$$
$$\mathcal{H}(z;n+1) - \mu(n) = (\mathcal{A}(z)\mathbf{E}_q + \mathcal{F}(z;n))(\mathcal{A}(z)\mathbf{E}_q + \mathcal{G}(z;n)) \quad (14)$$

*whith $\mathcal{F}$ and $\mathcal{G}$ such that the class of the Askey-Wilson polynomials becomes invariant under the action of the "raising" and "lowering" operators.*

Proof. The operatorial relations (14) are equivalent to the system:

$$\begin{aligned}
\mathcal{F}(qz;n) + \mathcal{G}(z;n) &= -(\mathcal{A}(qz) + \mathcal{B}(qz) + \mathcal{K}(qz)\lambda(n)) \\
\mathcal{F}(z;n)\mathcal{G}(z;n) &= \mathcal{A}(z)\mathcal{B}(qz) - \mu(n) \\
\Delta_q(\mathcal{F}(z;n) - \mathcal{G}(z;n)) &= (\lambda(n+1) - \lambda(n))\mathcal{K}(qz)
\end{aligned} \quad (15)$$

where $\Delta_q h(z) = h(qz) - h(z)$.
Using the $q$-integration, we first transform the system (15) in :

$$\begin{aligned}
\mathcal{F}(z;n) + \mathcal{F}(qz;n) &= -(\mathcal{A}(qz) + \mathcal{B}(qz)) - (\beta_{-1} - \tfrac{\lambda(n)}{q})z^{-1} - \beta_0 \\
&\quad - (\beta_1 + \lambda(n)q)z \\
\mathcal{F}(z;n)\mathcal{G}(z;n) &= \mathcal{A}(z)\mathcal{B}(qz) - \mu(n) \\
\mathcal{G}(z;n) - \mathcal{F}(z;n) &= \sum_{i=-1}^{i=1} \beta_i z^i
\end{aligned} \quad (16)$$

where $\beta_{-1} = \frac{\lambda(n+1)-\lambda(n)}{1-q}$; $\beta_1 = q\beta_{-1}$; $\beta_0$ remaining arbitrary for the moment.
Observing the first equation in (16) and then using the last one, it becomes sensible to search $\mathcal{F}(z;n)$ and then $\mathcal{G}(z;n)$ under the forms:

$$\begin{aligned}
\mathcal{F}(z;n) &:= \frac{F_{-2}z^{-2} + F_{-1}z^{-1} + F_0 + F_1 z + F_2 z^2}{qz - z^{-1}} \\
\mathcal{G}(z;n) &:= \frac{(F_{-2}-\beta_{-1})z^{-2} + (F_{-1}-\beta_0)z^{-1} + F_0 + (F_1+\beta_0 q)z + (F_2+\beta_1 q)z^2}{qz - z^{-1}}
\end{aligned} \quad (17)$$

Taking $\mu$, $\beta_0$, $F_{-2}$, $F_{-1}$, $F_0$, $F_1$, $F_2$ as unknowns, the system (15) will then be transformed in an algebraic system of 16 equations for 7 unknowns. To solve it (by hand), one first determines those unknowns from 7 equations and then very delicately ensures himself that they satisfy all the remaining 9 equations.
The result is:

$$\begin{aligned}
F_{-2}(n) &: & \frac{\lambda(n) - q\lambda(n+1)}{q^2 - 1} - \frac{q + A_2}{q^2 + q} \\
F_2(n) &: & \frac{\lambda(n)q - \lambda(n+1)}{1 - q^2}q^2 - \frac{q^2 + qA_2}{q+1}
\end{aligned} \quad (18)$$

$$\begin{aligned}
\beta_0(n) : &\quad \tfrac{1-q}{(\lambda(n)-\lambda(n+1))q^3}\{(2\tfrac{\lambda(n)q-\lambda(n+1)}{1-q^2}q^2 + \tfrac{\lambda(n+1)-\lambda(n)}{1-q}q^2 \\
&\quad -2\tfrac{q^2+qA_2}{1+q})(A_1+qA_{-1})+(2A_1q^2+2A_2A_{-1}q)\} \\
F_{-1}(n) : &\quad \tfrac{\beta_0(n)}{2} - \tfrac{A_1+qA_{-1}}{2q} \\
F_1(n) : &\quad -\tfrac{q\beta_0(n)}{2} - \tfrac{A_1+qA_{-1}}{2} \\
F_0(n) : &\quad \tfrac{1}{q+q^2}\{q^2-q^3-A_0(q+q^2)+A_2(q-1)+q^2(\lambda(n)+\lambda(n+1))\} \\
\mu(n) : &\quad A_0 + A_1A_{-1}q^{-1} + A_0A_2q^{-2} + F_0(n)\beta_{-1}(n) + F_{-1}(n)\beta_0(n) \\
&\quad -2F_{-2}(n)F_0(n) - F_{-1}^2(n).
\end{aligned}$$

Thus, (17) is determined and a solution of the IHF problem (14) is found. It remains to prove the invariance of the class of the AWP under the action of the "raising" and "lowering" operators thus obtained. From (14), the intertwining relations follow:

$$\begin{aligned}
\mathcal{H}(z;n+1)(\mathcal{A}(z)\mathbf{E}_q + \mathcal{F}(z;n)) &= (\mathcal{A}(z)\mathbf{E}_q + \mathcal{F}(z;n))\mathcal{H}(z;n) \\
\mathcal{H}(z;n)(\mathcal{A}(z)\mathbf{E}_q + \mathcal{G}(z;n)) &= (\mathcal{A}(z)\mathbf{E}_q + \mathcal{G}(z;n))\mathcal{H}(z;n+1)
\end{aligned} \quad (19)$$

Let us note that those relations, as well as the IHF (14) are in fact valid $\forall n \in Z$. So, if for a given number $n_0 \in Z$, a given function $\Psi_{n_0}(z)$ is an eigenvector of the operator $\mathcal{H}(z;n_0)$, corresponding to a given eigenvalue $\nu$ (take $\nu = 0$, in our situation), then the sequence of functions

$$^\uparrow\Psi_{n_0+k}(z) := \prod_{i=0}^{k-1}(\mathcal{A}(z)\mathbf{E}_q + \mathcal{F}(z;n_0+i))\Psi_{n_0}(z), k=1,2,3,\ldots \quad (20)$$

are respectively eigenvectors of $\mathcal{H}(z;n_0+k)$ corresponding to the same eigenvalue $\nu = 0$. The following formulas will be useful:

$$\mathcal{F}(z;n) - \mathcal{G}(z;n-1) =$$
$$[2A_2q^{-1}q^n - 2q^{-n}]\chi(z) - \tfrac{\beta_0(n)+\beta_0(n-1)}{2}, \quad (\forall n \in Z), \quad (21)$$
$$(\mathcal{A}(z)\mathbf{E}_q + \mathcal{F}(z;0))(1) = [2A_2q^{-1} - 2]\chi(z) - \tfrac{\beta_0(0)}{2} + \tfrac{A_1-qA_{-1}}{2q}, \quad (22)$$

Considering (20) and (14), we obtain the following difference relations (we let $^\uparrow\Psi_{n_0}(z) := \Psi_{n_0}(z)$):

$$\begin{aligned}
^\uparrow\Psi_{n+1}(z) &= (\mathcal{A}(z)\mathbf{E}_q + \mathcal{F}(z;n))^\uparrow\Psi_n(z) \\
-\mu(n)^\uparrow\Psi_n(z) &= (\mathcal{A}(z)\mathbf{E}_q + \mathcal{G}(z;n))^\uparrow\Psi_{n+1}(z) \\
& n = n_0, n_0+1, n_0+2, \ldots
\end{aligned} \quad (23)$$

Now, using (23) and (21), we obtain the following three-term recurrence relations:

$$^\uparrow\Psi_{n+1}(z) + \mu(n-1)^\uparrow\Psi_{n-1}(z) =$$
$$\left([2A_2q^{-1}q^n - 2q^{-n}]\chi(z) - \frac{\beta_0(n) + \beta_0(n-1)}{2}\right)^\uparrow \Psi_n(z) \quad (24)$$
$$n = n_0 + 1, n_0 + 2, \ldots$$

Considering the particular case of (20) when $\Psi_{n_0}(z)$ is the AWP of order $n_0$, $\mathcal{P}_{n_0}(\chi(z))$, the invariance of the class of the AWP under the "raising" operation means that $\Psi_{n_0+k}(z)$ will be a non vanishing constant multiple of $\mathcal{P}_{n_0+k}(\chi(z))$, the AWP of order $n_0 + k$, $k \geq 1$.

It follows from (22) and (24) (with $n_0 = 0$) that there exists a sequence of polynomial (in $\chi(z)$) eigenfunctions (of degrees $0, 1, 2, 3, \ldots$) of the Askey-Wilson operator $\mathcal{L}$, invariant under the "raising" operation. On the other side, according to (23), the same sequence is invariant, up to a multiplication by a constant, under the "lowering" operation (as long as $n \geq 0$). But as the constant multiple of AWP are the unique non trivial polynomial eigenfunctions of the Askey-Wilson operator, this suffices to show the cited invariance and the proposition is completely proved.

Thus, if in (20) we consider the case of AWP and set $n_0 = 0$, we obtain the Rodrigues type formula:

$$\tilde{c}(n)\mathcal{P}_n(\chi(z)) = \prod_{i=0}^{n-1} [\mathcal{A}(z)\mathbf{E}_q + \mathcal{F}(z;i)](1), n = 1, 2, 3, \ldots \quad (25)$$

$\tilde{c}(n)$ being some no null constant (in $z$). It is not difficult to see that (25) is equivalent to

$$\tilde{c}(n)\mathcal{P}_n(\chi(z)) = \frac{1}{\rho(z)} \prod_{i=0}^{n-1} [\mathbf{E}_q + \mathcal{F}(z;i)](\rho(z)), n = 1, 2, 3, \ldots \quad (26)$$

where

$$\frac{\rho(qz)}{\rho(z)} = \mathcal{A}(z),$$

a formula similar to that from [3].

Now, remembering that the relations (19) as well as the factorization (14) are valid $\forall n \in Z$, let us give the following notes for completness: Remark first that as in (20), the functions

$$^\downarrow\Psi_{n_0-1-k}(z) := \prod_{i=0}^{k}(\mathcal{A}(z)\mathbf{E}_q + \mathcal{G}(z; n_0 - 1 - i))\Psi_{n_0}(z), k = 0, 1, 2, \ldots \quad (27)$$

are respectively eigenfunctions of the operators $\mathcal{H}(z; n_0 - 1 - k)$. They satisfy the difference relations ( letting $^\downarrow\Psi_{n_0}(z) := \Psi_{n_0}(z)$):

$$\begin{aligned}
-\mu(n)^\downarrow\Psi_{n+1}(z) &= (A(z)\mathbf{E}_q + \mathcal{F}(z; n))^\downarrow\Psi_n(z) \\
^\downarrow\Psi_n(z) &= (A(z)\mathbf{E}_q + \mathcal{G}(z; n))^\downarrow\Psi_{n+1}(z) \qquad (28) \\
& n = n_0 - 1, n_0 - 2, \ldots
\end{aligned}$$

and the recurrence relations

$$\mu(n)^\downarrow\Psi_{n+1}(z) +^\downarrow\Psi_{n-1}(z) =$$
$$\left([2q^{-n} - 2A_2q^{-1}q^n]\chi(z) + \frac{\beta_0(n) + \beta_0(n-1)}{2}\right)^\downarrow\Psi_n(z) \quad (29)$$
$$n = n_0 - 1, n_0 - 2, \ldots$$

For $n_0 = 0$ and $\psi_0(z) = 1$ (polynomial case), (27) gives an "extension" of (25):

$$^\downarrow\Psi_{-(n+1)}(z) := \prod_{i=0}^{n}(\mathcal{A}(z)\mathbf{E}_q + \mathcal{G}(z; -(i+1)))(1), n = 0, 1, 2, \ldots \quad (30)$$

But a direct verification of the relation

$$(\mathcal{A}(z)\mathbf{E}_q + \mathcal{G}(z; -1))(1) = 0 \quad (31)$$

shows, without surprise, that the functions produced in (30), "extending" the AWP relatively to $n$ (as an index, from $Z^+$ to $Z^-$), are all null ($^\downarrow\Psi_{-(n+1)}(z) \equiv 0$, $n = 0, 1, 2, \ldots$).

The possibility of unifying the "raising" and "lowering" operations in (25) and (30) in types (23) or (28) pair of difference relations (and so in types (24) or (29) recurrence relations) with $n = \ldots - 2, -1, 0, 1, 2, \ldots$, is subordinated to the constraint $\mu(-1) = 0$ and this is actually the case.

We are now led to make a somewhat converse reasoning:
We remark first that

$$\begin{aligned}
\mathcal{F}(z;-1) &= -\mathcal{B}(qz) \\
\mathcal{G}(z;-1) &= -\mathcal{A}(z) \\
\mu(-1) &= 0
\end{aligned} \qquad (32)$$

This leads together with the first equation in (15) to a new form of the equation (12):

$$\{\mathcal{G}(qz;-1)\mathbf{E}_q^2 - [\mathcal{F}(qz;n) + \mathcal{G}(z;n)]\mathbf{E}_q + \mathcal{F}(z;-1)\}\mathcal{P}_n(\chi(z)) = 0. \qquad (33)$$

We recall that this equation is verified by $\mathcal{F}$ and $\mathcal{G}$ given by (18) and (17) and the AWP given by (25) or (26) for $n \geq 1$ while $\mathcal{P}_0 = 1$ (the unit is a solution of (33) for $n = 0$ thanks to the first equation in the system (5)). Let us return now the situation taking the IHF system (5) as a starting point. Our point resides in the following

**Proposition 2** *Consider the IHF system (5) with the unique initial condition $\tilde{\mu}(-1) = 0$ and let $\tilde{\mathcal{F}}$ and $\tilde{\mathcal{G}}$ be any of its solutions, then the equation*

$$\{\tilde{\mathcal{G}}(qz;-1)\mathbf{E}_q^2 - [\tilde{\mathcal{F}}(qz;n) + \tilde{\mathcal{G}}(z;n)]\mathbf{E}_q + \tilde{\mathcal{F}}(z;-1)\}\tilde{\mathcal{P}}_n(z) = 0 \qquad (34)$$

*admits a sequence of solutions given by a "Rodrigues type formula" (similar to (20)) applied to $\tilde{\mathcal{P}}_0 = 1$ and satisfying srict difference relations and three-term recurrence relations. Morever, if we add the condition*

$$\tilde{\mathcal{F}}(z;n) - \tilde{\mathcal{G}}(z;n-1) = c_0(n)\chi(z) + c_1(n) \qquad (35)$$

*for some constants (in z) $c_0$ and $c_1$, the obtained sequence is of polynomial (in $\chi(z)$) type.*

Proof. To prove this propostion, one needs essentially to remark that if $\tilde{\mathcal{F}}$ and $\tilde{\mathcal{G}}$ are solutions of the system (5) with initial condition $\tilde{\mu}(-1) = 0$, then the operator

$$\tilde{\mathcal{H}}(z;n) = \tilde{\mathcal{G}}(z;-1)\left[\tilde{\mathcal{G}}(qz;-1)\mathbf{E}_q^2 - [\tilde{\mathcal{F}}(qz;n) + \tilde{\mathcal{G}}(z;n)]\mathbf{E}_q + \tilde{\mathcal{F}}(z;-1)\right] \qquad (36)$$

admits a type (14) IHF. The remaining part of the reasoning is similar to that of the preceding proposition, reason for which we omit it and we assume that the proposition is proved.

It is clear that if we add the analogous of the two first initial conditions in (32), we obtain the AWP. In other words, we remark in our satisfaction that the IHF system (5) can be considered as an independent source of the AWP.

### 3. The generalization of the Askey-Wilson polynomials.

We envisage now to show how to produce a sequence of functions $\Phi_n(z,t)$ generalizing the AWP in that sense that they are eigenfunctions of the Askey-Wilson operator $\pounds$ (see (10)) and $\Phi_n(z,t) \longrightarrow \mathcal{P}_n(\chi(z))$ while $t \longrightarrow 1$. We showed above that if the eigenvalue of the Askey-Wilson operator is $\lambda(n) = -(1-q^{-n})(1-abcdq^{n-1})$, then an IHF as in (14) is realizable. Now, using (15), one finds that conversely, the possibility of a type (14) IHF implies necessarily that the eigenvalue in action is ${}^t\lambda(n) = -(1-tq^{-n})(1-abcdt^{-1}q^{n-1})$. On the other side, it is clear that the presence of the parameter $t$ in ${}^t\lambda(n)$ does not in no way annoy the IHF (14) (thanks to the specific dependence of ${}^t\lambda(n)$ at $t$). This means that if in the expressions (18), we substitute $\lambda(n)$ by ${}^t\lambda(n)$, the resulting expressions for $\mathcal{F}$ and $\mathcal{G}$, say ${}^t\mathcal{F}$ and ${}^t\mathcal{G}$, solve a type (14) IHF problem for an operator, say ${}^t\mathcal{H}(z;n)$, obtained from $\mathcal{H}(z;n)$ replacing $\lambda(n)$ by ${}^t\lambda(n)$.

Let us consider next the equation:

$$\pounds y(z,t) = {}^t\lambda(0) y(z,t), \tag{37}$$

or equivalently

$$ {}^t\mathcal{H}(z;0) y(z,t) = 0. \tag{38}$$

This equation is not new: A little wider equation (where the Askey-Wilson operator is replaced by the corresponding Nikiforov-Suslov-Uvarov operator [17]) has been explicitely solved in [2]. In our situation, according to [2], one first solves the non-homogeneous equations:

$$\left(\pounds - {}^t\lambda(0)\right) u(z,t,\alpha) = G(z,t,\alpha) \tag{39}$$

where

$$G(z,t,\alpha) = \frac{(\alpha t, abcd\alpha t^{-1}q^{-1};q)_\infty}{\alpha(ab\alpha, ac\alpha, ad\alpha, \alpha q;q)_\infty}(az, az^{-1};q)_\infty \tag{40}$$

and $(\sigma;q)_\infty := \lim_{i\to\infty}(\sigma;q)_i$, for the following values of $\alpha$:

$$\alpha := 1 \qquad \alpha := \frac{q}{ab} \qquad \alpha := \frac{q}{ac} \qquad \alpha := \frac{q}{ad}.$$

The corresponding solutions are the functions [2]:

$$u(z,t,\alpha) = \frac{(az, az^{-1}; q)_\infty}{(a\alpha z, a\alpha z^{-1}; q)_\infty} \sum_{i=0}^{\infty} \frac{(\alpha t, abcd\alpha t^{-1}q^{-1}, a\alpha z, a\alpha z^{-1}; q)_i}{(ab\alpha, ac\alpha, ad\alpha, \alpha q; q)_i} q^i. \quad (41)$$

Let us remark that the AWP appear here as the functions $u(z, q^{-n}, 1)$, $n = 0, 1, 2, \ldots$

The solutions of the homogeneous equation (37) or (38) are then obtained by opering adequate (considering (40)) linear combinations of any two of the functions (41).

Among the solutions of the homogenous equation (37) or (38), we are interested in those having a constant limit when $t$ converges to the unit. Let us take the "adequate" linear combination of $u(z, t, 1)$ and $u(z, t, \alpha)$, $\alpha \neq 1$:

$$y(z,t) := u(z,t,1) - \alpha \frac{(t, abcdt^{-1}q^{-1}, ab\alpha, ac\alpha, ad\alpha, \alpha q; q)_\infty}{(ab, ac, ad, q, \alpha t, abcd\alpha t^{-1}q^{-1}; q)_\infty} u(z,t,\alpha). \quad (42)$$

It is easily seen that $y(z, 1) = 1$. However, if $\tilde{y}(z, t)$ is another "adequate" linear combination of $u(z, t, \alpha_1)$ and $u(z, t, \alpha_2)$ where $\alpha_1$ and $\alpha_2$ differ both from the unit and each other, then $\tilde{y}(z, 1) \neq constant$. In other words $\tilde{y}(z, 1)$ is the nonconstant solution of the Askey-Wilson equation for $n = 0$.

Our point here resides in the following

**Proposition 3** *If $y(z, t)$ is the solution in (42), of (37) or (38), then the functions*

$$\begin{aligned}
\Phi_0(z,t) &:= y(z,t), \\
\Phi_n(z,t) &:= \prod_{i=0}^{n-1} [\mathcal{A}(z)\mathbf{E}_q + {}^t\mathcal{F}(z;i)] y(z,t), \\
&\quad n = 1, 2, 3, \ldots,
\end{aligned} \quad (43)$$

*generalize the AWP in that sense that they are eigenfunctions of the Askey-Wilson operator £ (corresponding to ${}^t\lambda(n)$) and they converge to them when $t$ converges to the unit.*

Proof. As already noted, the operator ${}^t\mathcal{H}$ admits a type (14) IHF but now with $\mathcal{F}$ and $\mathcal{G}$ replaced by ${}^t\mathcal{F}$ and ${}^t\mathcal{G}$. Consequentely, the fact that the functions in (43) are eigenfunctions of the Askey-Wilson operator £ corresponding to ${}^t\lambda(n)$, is a consequence (thanks to type (19) intertwining relations) of that $\Phi_0(z;t)$ is an eigenvector of £ corresponding to ${}^t\lambda(0)$. To be assured that the functions in (43) converge to the Askey-Wilson polynomials (of course, up to a multiplication by a constant), we need essentially to remember that ${}^t\lambda(n)$ converges to $\lambda(n)$ when $t$ converges to the unit and then compare the r.h.s. of (43) and (25), which proves the propositon.

For completness, It's perhaps worth noting that an iterative application of the "lowering" operator $\mathcal{A}(z)\mathbf{E}_q +{}^t \mathcal{G}(z;-(i+1))$, $i=0,1,2,\ldots$ to $y(z,t)$ leads to an evident "extension" relatively to $n$ (as an index, from $Z^+$ to $Z^-$) of the functions $\Phi_n(z,t)$ in (43), generalizing thus the already evoked identically null functions (30). Morever if ${}^\uparrow\tilde{\Psi}_n(z)$ (or ${}^\downarrow\tilde{\Psi}_n(z)$) are the eigenfunctions obtained by applying $n$-times the "raising" (or "lowering") operator to $\tilde{y}(z,1)$, so their generalizations are obtained by performing similar operations starting now at $\tilde{y}(z,t)$. Let us note finally that the present generalizations are also difference hypergeometric functions (in the sense of [2]) reason for which they are closely related to the AWP.

**Conclusion.** We have (mainly) shown explicitly in this work and in [3] that all the classical (up to the AWP) orthogonal polynomials can be recovered solving the IHF system (5) ($n$ being presumed to be the degree of the searched polynomials). We have serious indications allowing to believe that there (i.e. a little modification of the IHF system (5)) is exactely the central "front door" when attempting to handle concretely non-classical classes of orthogonal polynomials such as the discrete semi-classical orthogonal polynomials or even the discrete Laguerre-Hahn orthogonal polynomials (the largest nowadays known class of orthogonal polynomials) [11], [12], [13]. We reserve a concrete discussion in this matter for a later publication.

**Aknowledgments.**

We would like to thank Professor A P Magnus for suggesting the present explorations and for fruitful discussions. Thanks are addressed as well to the Belgian General Agency for Cooperation with Developing Countries (AGCD) for financial support.

# References

[1] R. Askey, J. Wilson, Some basic hypergeometric orthogonal polynomials that generalize Jacobi polynomials, *Mem. Am. Math. Soc.* **54** (1985)

1-55.

[2] N. M. Atakishiyev, S. K. Suslov, Difference hypergeometric functions, in: A. A. Gonchar and E. B. Saff, Eds., *Prog. Approx. Theory* (Springer-Verlag, 1992) 1-35.

[3] G. Bangerezako, Discrete Darboux transformation for discrete polynomials of hypergeometric type , *J. Phys. A: Math. Gen.* **31** (1998) 1-6.

[4] B. Brown, W. Evans, M. Ismail, The Askey-Wilson polynomials and $q$-Sturm-Liouville problems, *Math. Proc. Camb. Phil. Soc.* **119** No1 (1996) 1-16.

[5] G. Darboux, *Théorie des surfaces II, p. 213* ( Gauthier-Villars, Paris 1915).

[6] B. Grammaticos, A. Ramani, Discrete Painlevé equations: derivation and proprieties, *NATO Adv. Sci. Inst. Ser. C Math. Phys. Sci.* **413** (1993) 299-313.

[7] L. Infeld, T. Hull, The factorization method, *Rev. Mod. Phys.* **23** (1951) 21-68.

[8] T. Inui, Unified theory of recurrence formulas I, *Prog. Theor. Phys.* **3** No 2 (1948) Apr-June.

[9] E. G. Kalnins and W. Miller, Symmetry techniques for $q$-series: Askey-Wilson polynomials, *Rocky Mountain J. Math.* **19** (1989) 223-230.

[10] R. Koekoek, R. Swarttouw, The Askey-scheme of hypergeometric orthogonals and its $q$-analogue, *Report of the Technical University Delft, Faculty of Technical Mathematics and Informatics* (1994) 94-05.

[11] A. P. Magnus, Riccati acceleration of Jacobi continued fractions and Laguerre-Hanh orthogonal polynomials, *Lect. Notes. Math.* **1071** (Springer, Berlin, 1984) 213-230.

[12] A. P. Magnus, Associated Askey-Wilson polynomials as Laguerre-Hahn orthogonal polynomials, *Springer Lect. Notes in Math.* **1329** (Springer, Berlin, 1988) 261-278.

[13] A. P. Magnus, Special non uniform lattice (snul) orthogonal polynomials on discrete dense sets of points, *J. Comp. Appl. Math.* **65** (1997) 253-265.

[14] V. Matveev, M. Salle, Differential-difference evolution equations II, Darboux transformation for the toda lattice, *Lett. Math. Phys.* **3** (1979) 425-429.

[15] W. Miller, Lie theory and difference equations I, *J. Math. Anal. Appl.* **28** (1969) 383-99.

[16] W. Miller, Symmetry techniques and orthogonality for $q$-series, *IMA volumes in Math. Appl.* **18** (1989) 191-212.

[17] A. Nikiforov, S. Suslov, V. Uvarov, *Classical orthogonal polynomials of a discrete variable* (Spring-Verlag, Berlin-Heidelberg-New York, 1991).

[18] A. B. Shabat, A. P. Veselov, Dressing chain and spectral theory of Schrdinger operator, *Funct. Anal. Appl.* **27** (1993) 81-96.

[19] Yu. F. Smirnov, On factorization and algebraization of difference equations of hypergeometric type, in: M. Alfaro and al., Eds., *Proceedings of the International Workshop on Orthogonal Polynomials in Mathematical Physics (in honour of Professor André Ronveaux)* (Legans june 24-26, 1996).

[20] V. Spiridonov, Symmetries of factorization chains for the discrete Shrödinger equation, *J. Phys. A: Math. Gen.* **30** (1997) L15-L21.

[21] V. Spiridonov, L. Vinet, A. Zhedanov, Difference Schrdinger operators with linear and exponential discrete spectra, *lett. Math. Phys.* **29** (1993) 63-73.

[22] V. Spiridonov, L. Vinet, A. Zhedanov, Periodic reductions of the factorization chain and the Hahn polynomials, *J. Phys. A: Math. Gen.* **27** (1994) L669-675.

[23] V. Spiridonov, L. Vinet, A. Zhedanov, Discrete Darboux transformations, the discrete-time Toda lattice and the Askey-Wilson polynomials, *Methods Appl. Anal.* **2** (1995) 369-398.

Gaspard Bangerezako
Institut de Mathmatique, Universit Catholique de Louvain,
Chemin du Cyclotron 2, B-1348 Louvain-La-Neuve, Belgium
*fax*: 32 10 47 25 30, *e-mail* :bangerezako@agel.ucl.ac.be